\documentclass{article}

\usepackage{amssymb}
\usepackage{amsmath,amsthm}
\usepackage[dvips]{graphicx}
\usepackage{wrapfig}
\usepackage{bm}

\newtheorem{theorem}{Theorem}[section]
\newtheorem{corollary}[theorem]{Corollary}
\newtheorem{lemma}[theorem]{Lemma}

\theoremstyle{definition}

\theoremstyle{remark}

\begin{document}
\title{Lower bounds for the warping degree of a knot projection}

\author{Atsushi Ohya\thanks{Department of Computer Science and Engineering, University of Yamanashi, 4-4-37, Takeda, Kofu-shi, Yamanashi, 400-8510, Japan. }
\and 
Ayaka Shimizu\thanks{Department of Mathematics, National Institute of Technology (KOSEN), Gunma College, 580 Toriba, Maebashi-shi, Gunma, 371-8530, Japan. Email: shimizu@gunma-ct.ac.jp, shimizu1984@gmail.com}}

\maketitle

\begin{abstract}
The warping degree of an oriented knot diagram is the minimal number of crossings which we meet as an under-crossing first when we travel along the diagram from a fixed point. 
The warping degree of a knot projection is the minimal value of the warping degree for all oriented alternating diagrams obtained from the knot projection. 
In this paper, we consider the maximal number of regions which share no crossings for a knot projection with a fixed crossing, and give lower bounds for the warping degree.  
\end{abstract}

\section{Introduction}

\begin{figure}[ht]
\begin{center}
\includegraphics[width=20mm]{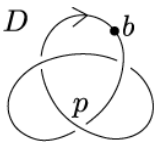}
\caption{The oriented based knot diagram $D_b$ has warping degree one. The crossing $p$ is the warping crossing point of $D_b$. }
\label{wcp}
\end{center}
\end{figure}

In this paper we assume that every knot diagram and knot projection has at least one crossing. 
A {\it based knot diagram} is a knot diagram which is given a base point on the diagram avoiding crossings. 
We denote by $D_b$ a based diagram $D$ with the base point $b$. 
The {\it warping degree}, $d(D_b)$, of an oriented based knot diagram $D_b$ is the number of crossings such that we encounter the crossing as an under-crossing first when we travel along $D$ with the orientation starting at $b$. 
We call such a crossing a {\it warping crossing point of $D_b$} (see Figure \ref{wcp}). 
The {\it warping degree}, $d(D)$, of an oriented knot diagram $D$ is the minimal value of $d(D_b)$ for all base points $b$ of $D$ (\cite{kawauchi-lecture}). 
A knot diagram is said to be {\it monotone}, or descending, if the warping degree is zero. 
Conversely, we can assume that the warping degree represents a complexity of a diagram in terms of how distant a knot diagram is from a monotone diagram. 
Note that a monotone knot diagram is a diagram of the trivial knot. 
A knot diagram is said to be {\it alternating} if we encounter an over-crossing and an under-crossing alternatively when traveling the diagram starting at any point on the diagram.

\begin{figure}[ht]
\begin{center}
\includegraphics[width=120mm]{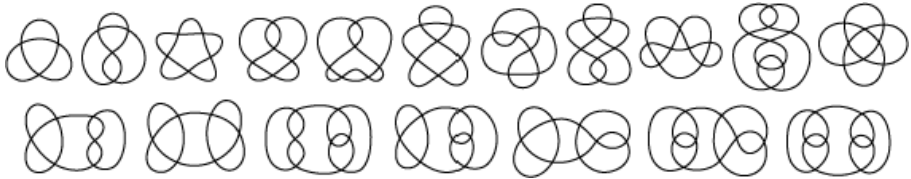}
\caption{All the reduced knot projections of warping degree one or two. The first two knot projections have warping degree one, and the others have two. }
\label{wd12}
\end{center}
\end{figure}

Let $P$ be an unoriented knot projection. 
The {\it warping degree of $P$} is defined to be the minimal value of the warping degree for all the oriented alternating diagrams obtained from $P$ by giving the orientation and crossing information. 
In Figure \ref{wd12}, all the reduced knot projections with warping degree one and two are shown (\cite{shimizu-wd-p}). 
Further examples for warping degree three or four are listed in the table in Section \ref{section-table} in this paper. 
As we may see, the warping degree of a knot projection shows somewhat complexity of a knot projection, like how ``curly'' a knot projection is, or how ``quick'' to back to a crossing when traveling the projection. 

The knot projections of warping degree two are determined in \cite{shimizu-wd-p} by considering all possibilities of connections of the unavoidable parts. 
Further explorations in the same way for warping degree three or more would be difficult since there are too many kinds of unavoidable parts and too many possibilities of their connections. 
In this paper, we introduce the {\it maximal independent region number}, $\mathrm{IR}(P)$, of a knot projection in Section \ref{section-IRS}, and show the following inequality which is useful to estimate the warping degree. 

\phantom{x}
\begin{theorem}
The inequality
$$\mathrm{IR}(P) \leq d(P) \leq c(P)-\mathrm{IR}(P)-1$$
holds for every reduced knot projection $P$, where $c(P)$ denotes the crossing number of $P$. 
\label{IRdIR}
\end{theorem}
\phantom{x}

\noindent As mentioned in Sections \ref{section-IRS} and \ref{section-IandR}, the value of $\mathrm{IR}(P)$ can be obtained without traveling along the knot projection, and also calculated just by solving simultaneous equations. \\

\begin{figure}[ht]
\begin{center}
\includegraphics[width=60mm]{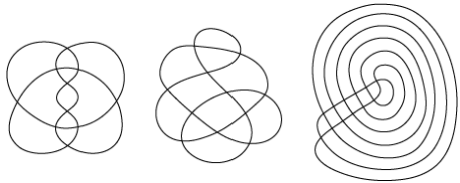}
\caption{Knot projections of warping degree three with 10, 11, 12 crossings. }
\label{three}
\end{center}
\end{figure}
\noindent Since all the reduced knot projections with warping degree one and two are determined and we can find some knot projections with warping degree three (see Figure \ref{three} and Section \ref{section-table}), we obtain the following table about the minimal value of the warping degree of reduced knot projections for each crossing number. 
\begin{table}[ht]
\begin{center}
\begin{tabular}{|c||c|c|c|c|c|c|c|c|c|c|} \hline
$c$ & 3 & 4 & 5 & 6 & 7 & 8 & 9 & 10 & 11 & 12 \\ \hline
$d^{\text{min}}(c)$ & 1 & 1 & 2 & 2 & 2 & 2 & 3 & 3 & 3 & 3 \\ \hline
\end{tabular}
\caption{The crossing number $c$ and the minimal value of warping degree $d^{\text{min}}(c)$ for all reduced knot projections with $c$ crossings. }
\end{center}
\end{table}
Regions of a knot or link projection are {\it independent} if they share no crossings. 
We also give the following lower bound for the warping degree which would be helpful to extend the above table. 

\phantom{x}
\begin{theorem}
If all the connected link projections with $n$, $n+1$ or $n+2$ crossings have $m$ or more independent regions, then $d(P) \geq m-1$ holds for all reduced knot projections $P$ with $n$ or more crossings. 
\label{thm-m-1}
\end{theorem}
\phantom{x}

\noindent The rest of the paper is organized as follows: 
In Section \ref{section-IRS}, we define the maximal independent region number $\mathrm{IR}(P)$ and prove Theorem \ref{IRdIR}. 
In Section \ref{section-IandR}, we introduce the calculation for $\mathrm{IR}(P)$ by simultaneous equations. 
In Section \ref{section-EM}, we estimate $\mathrm{IR}(P)$ and the warping degree $d(P)$ and prove Theorem \ref{thm-m-1}. 
In Section \ref{section-table}, we list and compare the values of $\mathrm{IR}(P)$ and $d(P)$.

\section{Independent region sets}
\label{section-IRS}

In this section, we define the maximal independent region number, and using it we estimate the warping degree of a knot projection. 
Throughout this section, we assume that every knot diagram and knot projection is reduced. 
We have the following lemma (cf. \cite{shimizu-wd-p}). 

\phantom{x}
\begin{lemma}
Let $D$ be an oriented alternating knot diagram. 
Let $c$ be a crossing of $D$. 
Take a base point $b$ just before an over-crossing of $c$. 
If $D$ has a region $R$ which does not incident to $c$, then one of the crossings on the boundary of $R$ is a warping crossing point of $D_b$, and one of that is a non-warping crossing point of $D_b$. 
\label{lem-bcR}
\end{lemma}
\phantom{x}

\begin{proof}
Let $e$ be the edge on the boundary of $R$ such that we meet it first from $b$. 
Then $e$ has an under-crossing and an over-crossing, that is, a warping crossing point and a non-warping crossing point. 
\begin{figure}[ht]
\begin{center}
\includegraphics[width=50mm]{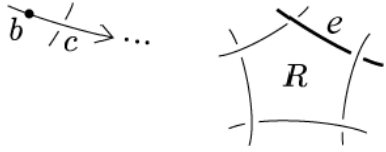}
\caption{The edge $e$ on the boundary of $R$ which we meet first from the base point $b$ has one under-crossing and one over-crossing, and they are a warping crossing point and non-warping crossing point of $D_b$, respectively, regardless of the orientation.  }
\label{bcR}
\end{center}
\end{figure}
\end{proof}
\phantom{x}

\noindent Similarly, we have the following. 

\phantom{x}
\begin{corollary}
Let $D$ be an oriented alternating knot diagram. 
Let $c$ be a crossing of $D$. 
Take a base point $b$ just before an over-crossing of $c$. 
If $D$ has independent $n$ regions which are not incident to $c$, then the inequality $n \leq d(D) \leq c(D)-n-1$ holds. 
\end{corollary}
\phantom{x}

\begin{proof}
By Lemma \ref{lem-bcR}, $D$ has at least $n$ warping crossing points of $D_b$. 
Also, $D$ has at least $n+1$ non-warping crossing points since the crossing $c$ is a non-warping crossing point, too. 
Therefore we have $n \leq d(D_b) \leq c(D)-n-1$. 
By the location of the base point $b$, we have $d(D_b)=d(D)$ (\cite{shimizu-wd}). 
\end{proof}
\phantom{x}

\noindent By definition, we have the following corollary for knot projections. 

\phantom{x}
\begin{corollary}
Let $P$ be a knot projection, and $c$ a crossing of $P$. 
If $P$ has independent $n$ regions which are not incident to $c$, then the inequality $n \leq d(P) \leq c(P)-n-1$ holds. 
\label{cor-ndc}
\end{corollary}
\phantom{x}

\begin{figure}[ht]
\begin{center}
\includegraphics[width=60mm]{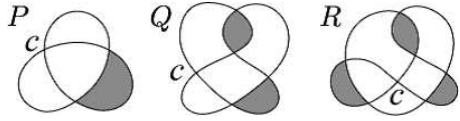}
\caption{The warping degree can be determined for some knot projections from the inequality of Corollary \ref{cor-ndc}, without traveling along the knot projection. The knot projection $P$ has $1 \leq d(P) \leq 3-1-1$, and $d(P)=1$. We also obtain $d(Q)=2, d(R)=3$ from the inequalities. }
\label{IRd-ex}
\end{center}
\end{figure}

\noindent The strong point is that we can estimate the warping degree {\it without traveling the projection} (see Figure \ref{IRd-ex}). 
This would enable us to estimate the warping degree more combinatorically. 
We call the set of regions of a knot projection $P$ which are independent and are not incident to a crossing $c$ an {\it independent region set for $P^c$}. 
We call the crossing $c$ a {\it base crossing}. 
We define the {\it maximal independent region number of $P^c$}, $\mathrm{IR}(P^c)$, to be the maximal cardinality of an independent region set for $P^c$. 
We define the {\it maximal independent region number of $P$}, $\mathrm{IR}(P)$, to be the maximal value of $\mathrm{IR}(P^c)$ for all base crossings $c$. 
We prove Theorem \ref{IRdIR}. \\

\noindent {\it Proof of Theorem \ref{IRdIR}.} \ 
It follows from Corollary \ref{cor-ndc}. 
\hfill$\Box$

\phantom{x}

\noindent

\section{Independent region sets and region choice matrix}
\label{section-IandR}

\begin{figure}[ht]
\begin{center}
\includegraphics[width=70mm]{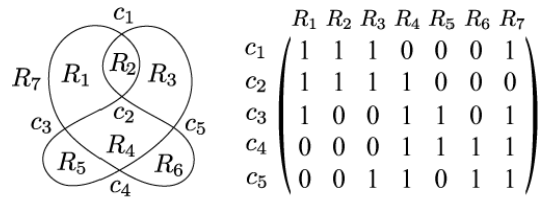}
\caption{A region choice matrix $M=(m_{i j})$, where $m_{i j}=1$ if $R_j$ is incident to $c_i$, and otherwise $m_{i j}=0$. }
\label{rcmatrix}
\end{center}
\end{figure}

In this section we explore how to find the independent region sets. 
A {\it region choice matrix} $M$, defined in \cite{AS}, of a knot projection $P$ of $n$ crossings is the following $n \times (n+2)$ matrix. 
(The transposition is known as an {\it incidence matrix} defined in \cite{CG}.) 
If a crossing $c_i$ is on the boundary of a region $R_j$, the $(i,j)$ component of $M$ is 1, and otherwise 0 (see Figure \ref{rcmatrix}). 

Now we find out all of the independent region sets for $P^{c_3}$ for the knot projection $P$ and the crossing $c_3$ in Figure \ref{rcmatrix} by looking at its region choice matrix. 
Since the crossing $c_3$ is involved with the four regions $R_1, R_4, R_5$ and $R_7$, namely, the third row has 1 at the first, fourth, fifth and seventh column, we can not choose them as independent regions for $P^{c_3}$. 
Hence we choose the regions from the rest regions $R_2, R_3$ and $R_6$. 
Namely, we choose columns from the second, third and sixth so that there are no components with the value two or more in the sum of the columns. 
Thus we obtain all the independent region sets for $P^{c_3}$, as $\{ R_2 \}$, $\{ R_3 \}$, $\{ R_6 \}$ and $\{ R_2, R_6 \}$. 

More generally, we can find out all the independent region sets for a knot projection $P^c$ for all base crossings $c$ from the region choice matrix by solving the following simultaneous equations\footnote{
If it works on $\mathbb{Z}_2$, it is known that the simultaneous equations have solutions for any $b_i$'s and any region choice matrix of a knot projection (\cite{shimizu-rcc}, \cite{CG}). 
Besides, if $x_i$'s are permitted to have the value for any integer, it is also known that the simultaneous equations have solutions for any region choice matrix of a knot projection even if $b_i$'s have the value for any integers (\cite{AS}). 
In this case, however, the equation has no solutions for some $b_i$'s. 
In Lemma \ref{lem-geq1}, we will see that it definitely has solutions for some $b_i$'s.}
 for $x_i \in \{ 1,0 \} \ (i=1,2, \dots ,7)$ 

\begin{align*}
\left(
\begin{array}{ccccccc}
1 & 1 & 1 & 0 & 0 & 0 & 1 \\
1 & 1 & 1 & 1 & 0 & 0 & 0 \\
1 & 0 & 0 & 1 & 1 & 0 & 1 \\
0 & 0 & 0 & 1 & 1 & 1 & 1 \\
0 & 0 & 1 & 1 & 0 & 1 & 1 
\end{array}
\right)  \left(
\begin{array}{c}
x_1 \\
x_2 \\
x_3 \\
x_4 \\
x_5 \\
x_6 \\
x_7
\end{array}
\right) = \left(
\begin{array}{c}
b_1 \\
b_2 \\
b_3 \\
b_4 \\
b_5 
\end{array}
\right) ,
\end{align*}

\noindent for all $b_i \in \{ 1,0 \} \ (i=1,2, \dots ,7)$, where $b_k \ne b_l$ for some $k$ and $l$; 
If $b_i=0$ for all $i$, it implies that no regions are chosen. 
If $b_i=1$ for all $i$, it means all the crossings are on the chosen regions, and we can not have a base crossing.

\section{Estimation for the maximal independent region number}
\label{section-EM}

As Theorem \ref{IRdIR} implies, the warping degree is estimated by the maximal independent region numbers. 
In this section, we estimate the maximal independent region number itself, and prove Theorem \ref{thm-m-1}. 
First, we have the following: 

\phantom{x}
\begin{lemma}
The inequality
\begin{align*}
\mathrm{IR}(P) \leq \frac{c(P)-1}{2}
\end{align*}
holds for every reduced knot projection $P$. 
\end{lemma}
\phantom{x}

\begin{proof}
By Theorem \ref{IRdIR}, we have $\mathrm{IR}(P) \leq c(P)-\mathrm{IR}(P) -1$, and then have $2\mathrm{IR}(P) \leq c(P)-1$. 
\end{proof}
\phantom{x}

\noindent Next, we give a routine lower bound for $\mathrm{IR}(P)$. 

\phantom{x}
\begin{lemma}
For any knot projection $P$ with $c(P) \geq 2$, we have $\mathrm{IR}(P) \geq 1$. 
\label{lem-geq1}
\end{lemma}
\phantom{x}

\begin{proof}
For the case that $c(P)=2$, $P$ has two independent bigons. 
By taking a base crossing at a crossing which belongs to one of the two bigons, we can take an independent region at the other bigon. 
For the case that $c(P) \geq 3$, then the number of regions is $3+2=5$ or more by the Euler characteristic (see, for example, \cite{AS}). 
Take a base crossing $c$. 
Then either three or four regions are incident to $c$. 
This means there exists a region which is not incident to $c$. 
Hence $\mathrm{IR}(P) \geq 1$ holds. 
\end{proof}
\phantom{x}

\noindent To give further lower bounds for $\mathrm{IR}(P)$, we show the following lemma for link projections. 

\phantom{x}
\begin{lemma}
If all the connected link projections with $n$, $n+1$ or $n+2$ crossings have $m$ or more independent regions, then all the connected link projections with $n$ or more crossings have $m$ or more independent regions. 
\label{lem-n+2}
\end{lemma}
\phantom{x}

\begin{proof}
Let $P$ be a link projection with $n+3$ crossings. 
If $P$ is reducible, splice it at a reducible crossing as shown in Figure \ref{reducible}. 
\begin{figure}[ht]
\begin{center}
\includegraphics[width=85mm]{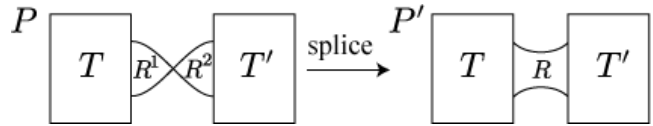}
\caption{Splice $P$ at a reducible crossing. }
\label{reducible}
\end{center}
\end{figure}
Then we obtain a link projection, $P'$, with $n+2$ crossings. 
By assumption, $P'$ has $m$ independent regions. 
Take the corresponding regions of $P$; 
If the region $R$ of $P'$ created by the splice has been chosen, take one of the parts $R^1$ and $R^2$ as a corresponding region (see Figure \ref{reducible}). 
Thus, we obtain $m$ independent regions of $P$. \\

If $P$ is reduced, it is shown in \cite{AST} that $P$ has a bigon or trigon. 
Splice it at a bigon or a trigon as shown in Figure \ref{reduced}. 
\begin{figure}[ht]
\begin{center}
\includegraphics[width=90mm]{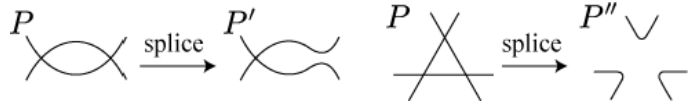}
\caption{Splice $P$ at a bigon or a trigon and obtain a link projection $P'$ or $P''$, respectively. Note that $P'$ has $n+2$ and $P''$ has $n$ crossings. }
\label{reduced}
\end{center}
\end{figure}
Then $P'$ and $P''$ have $m$ independent regions. 
Similarly, take $m$ regions of $P$ properly from the corresponding regions, which are independent. 
\end{proof}
\phantom{x}

\noindent We have the following corollary. 

\phantom{x}
\begin{corollary}
If all the connected link projections with $n$, $n+1$ or $n+2$ crossings have $m$ or more independent regions, then $\mathrm{IR}(P) \geq m-1$ holds for all reduced knot projections $P$ with $n$ or more crossings. 
\label{cor-knotm}
\end{corollary}

\phantom{x}
\begin{proof}
All reduced knot projections with $n$ or more crossings have $m$ independent regions
by Lemma \ref{lem-n+2}. 
Take a base crossing $c$ at the boundary of one of the $m$ regions. 
Then, the rest $m-1$ regions are independent regions for $P^c$. 
\end{proof}
\phantom{x}

\noindent We prove Theorem \ref{thm-m-1}. \\

\noindent {\it Proof of Theorem \ref{thm-m-1}. } \ 
From Corollary \ref{cor-knotm} we have $\mathrm{IR}(P) \geq m-1$, and from Theorem \ref{IRdIR} we have $d(P) \geq \mathrm{IR}(P)$. 
\hfill$\Box$

\section{Table of $\mathrm{IR}(P)$ and $d(P)$}
\label{section-table}

For the knot projections $P$ of prime alternating knots with up to nine crossings which are obtained from the knot diagrams in Rolfsen's knot table (\cite{rolfsen}), the values of $\mathrm{IR}(P)$ and $d(P)$ (\cite{shimizu-wd}) are listed below. 
The values of $\mathrm{IR}(P)$ were obtained by the calculation using the SAT solver. 

\begin{figure}[ht]
\begin{center}
\includegraphics[width=120mm]{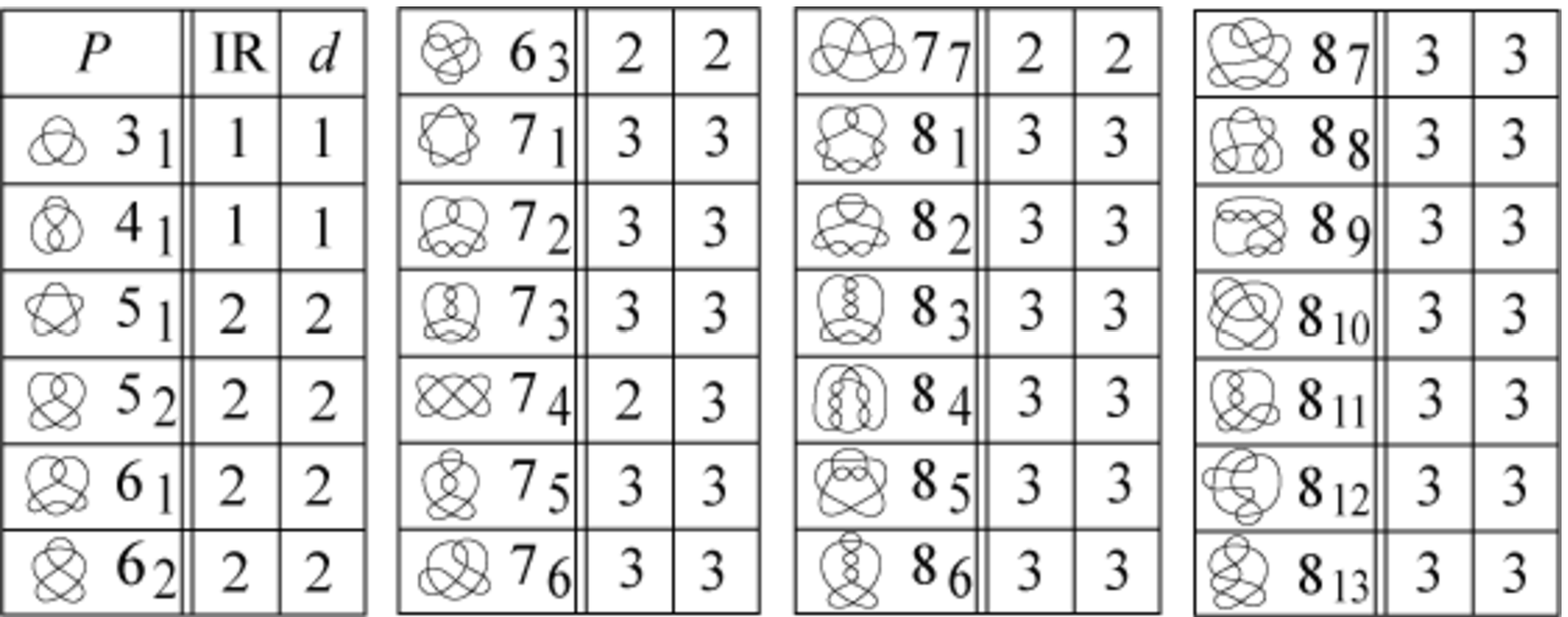}
\label{table1}
\end{center}
\end{figure}

\begin{figure}[ht]
\begin{center}
\includegraphics[width=120mm]{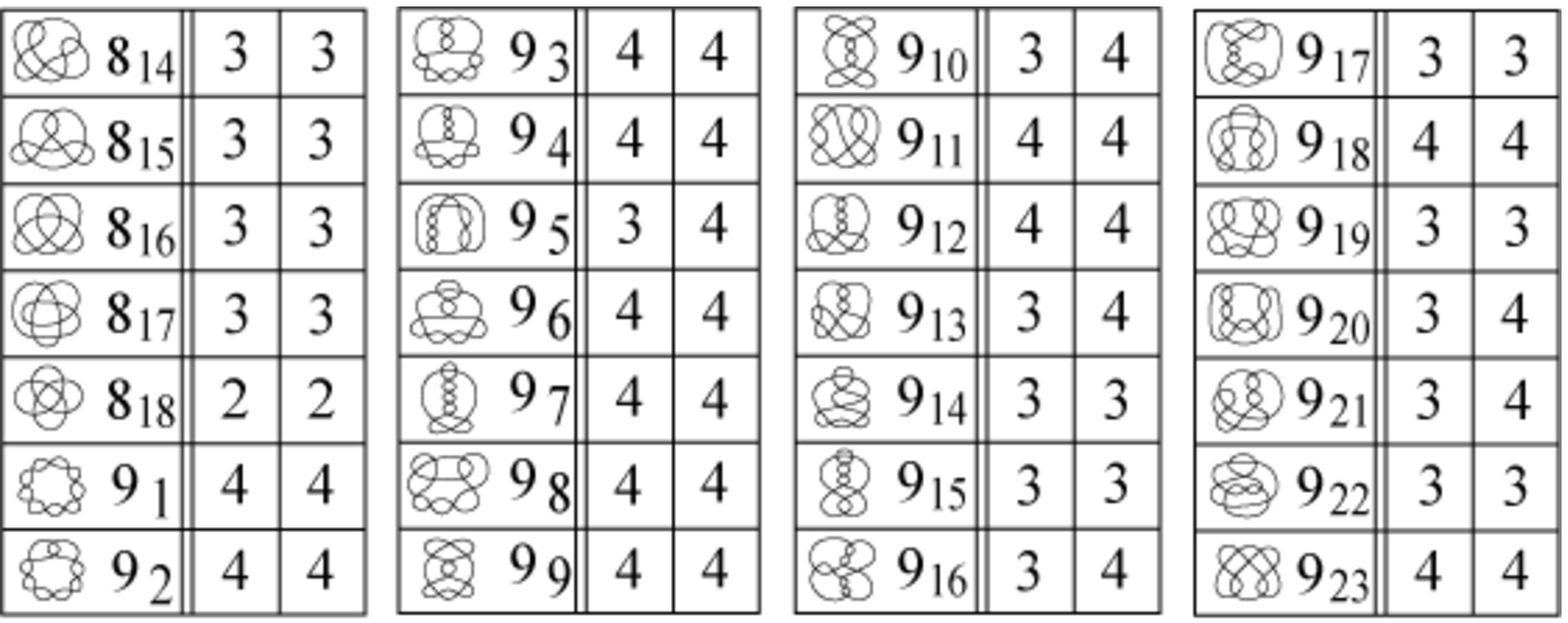}
\label{table2}
\end{center}
\end{figure}

\begin{figure}[ht]
\begin{center}
\includegraphics[width=86mm]{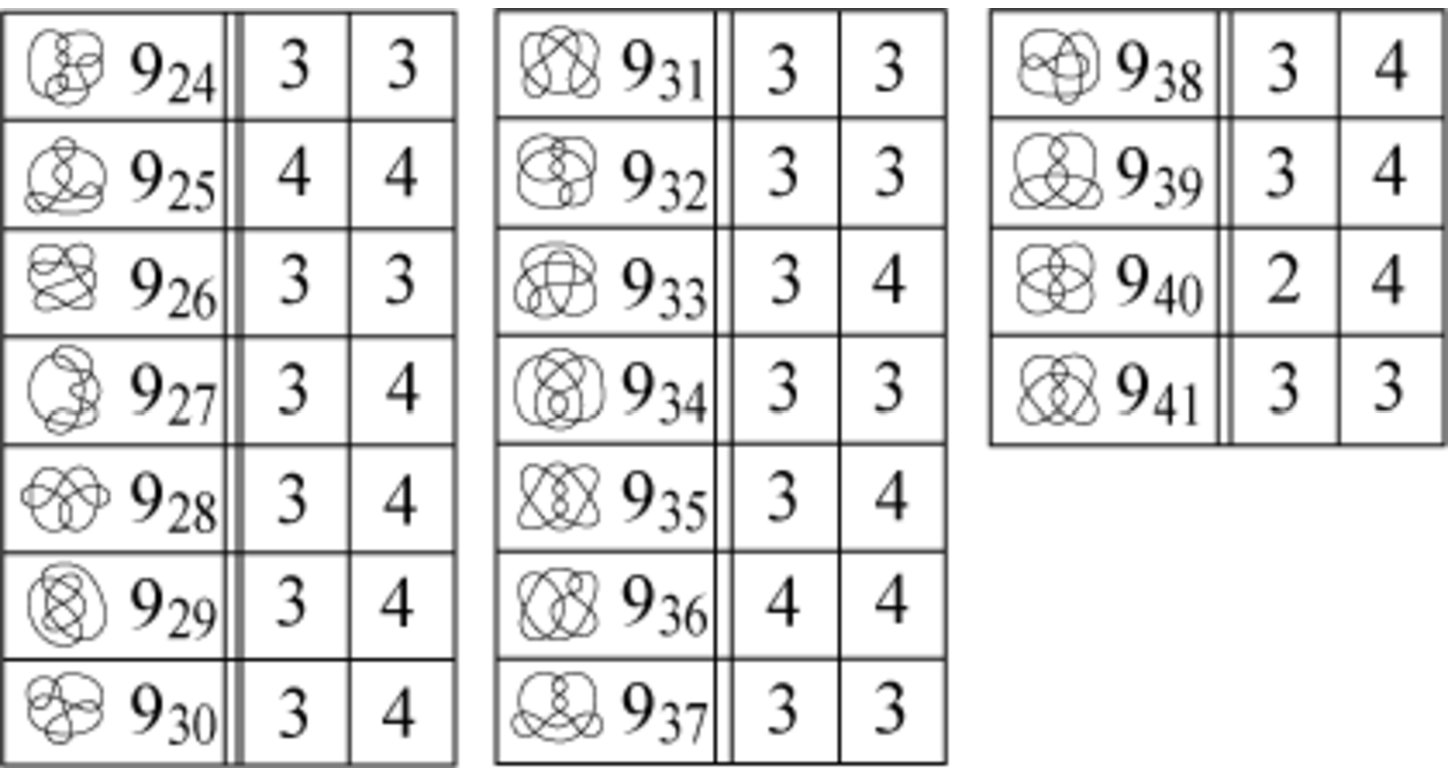}
\label{table3}
\end{center}
\end{figure}

\section*{Acknowledgment}
The authors thank Yoshiro Yaguchi for helpful comments.

\end{document}